\magnification=\magstep1 
\input amssym.def
\baselineskip=13pt
\def \ab {^{\rm ab}}
\def \s {^{\rm s}}
\def \E {{\cal E}}
\def \Z {{\Bbb Z}}
\def \F {{\Bbb F}}
\def \R {{\Bbb R}}
\def \Q {{\Bbb Q}}
\def \C {{\Bbb C}}
\def \P {{\Bbb P}}
\def \A {{\Bbb A}}
\def \G {{\Bbb G}}
\def \Gal {{\rm Gal}}
\def \ker {{\rm ker}}
\def \Hom {{\rm Hom}}
\def \iso {\to^{\!\!\!\!\!\!\!\sim\,}}
\def \card {{\rm card}\,}
\def \cd {{\rm cd}}
\def \et {{\hbox{\sevenrm \'et}}}
\def \Spec {{\rm Spec}}
\def \Br {{\rm Br}}
\def \div {{\rm div}}
\def \red {{\rm red}}
\def \char {{\rm char}}
\def\sdp{{\times\kern-.27em\raise1pt
  \hbox{$\scriptscriptstyle{|}$}
  {\kern-.25em\raise1pt\hbox{$\scriptscriptstyle{|}$}}
  \kern.1666em}}
\def \qed {\hfill\lower0.9pt\vbox{\hrule \hbox{\vrule height 0.2 cm   
  \hskip 0.2 cm \vrule height 0.2 cm}\hrule}} 
\font\small=cmr8
\font\smit=cmti8

\centerline{\bf On Function Fields with Free Absolute Galois Groups} 
\vskip .4cm
\centerline{David Harbater\footnote{$^*$}{\small Supported in part by
NSF Grant DMS-0500118.
{\parindent=0pt\item{}  
{\smit 2000 Mathematics Subject Classification}.  Primary 12E30, 12F10, 14H30; Secondary 12D15, 14H05, 20E05.  
\item{}  \baselineskip=10pt
{\smit Key words and phrases}: absolute Galois group, embedding problem, Galois covers, profinite groups, function fields, quasi-free group.}}}
\centerline{Dept.\ of Mathematics, University of Pennsylvania}

\bigskip

{\narrower\noindent{\bf Abstract.} We prove that certain fields have the property that their absolute Galois groups are free as profinite groups: the function field of a real curve with no real points; the maximal abelian extension of a 2-variable Laurent series field over a separably closed field; and the maximal abelian extension of the function field of a curve over a finite field.  These results are related to generalizations of Shafarevich's conjecture. \par}

\baselineskip=13pt

\medskip
 
\noindent  {\bf  \S1: Introduction.}

\medskip

This paper shows that the absolute Galois groups of certain fields are free as profinite groups.  Although these fields arise in geometric contexts, our results are related to Shafarevich's conjecture on absolute Galois groups, which he posed in the context of number theory.  In its original form, that conjecture states that the absolute Galois group of $\Q\ab$ is free, where $\Q\ab$ is the maximal abelian extension of $\Q$.  Since $\Q\ab$ is also the maximal cyclotomic extension of $\Q$, Shafarevich's conjecture has been generalized to assert that for any global field $K$, the absolute Galois group of the maximal cyclotomic extension of $K$ is free (see e.g.\ [Ha2]).  That conjecture remains open in the number field case (even for $\Q$), but it was proven in the function field case in [Ha1] and [Po1].  In that case, the conjecture is equivalent to saying that for any curve over $\bar\F_p$, the absolute Galois group of the function field is free.  In [Ha1] and [Po1] even more was shown: that freeness holds for any curve over {\it any} algebraically closed field.  

The result in [Ha1] and [Po1] suggests asking what happens over fields $K$ that are ``almost'' algebraically closed; i.e.\ such that $[\bar K:K]$ is finite, where $\bar K$ is the algebraic closure.  By the theorem of Artin-Schreier, these are precisely the real closed fields (e.g.\ $\R$); and in Theorem~4.2 we show that the function field of a curve $X$ over a real closed field $R$ has free absolute Galois group if and only if $X$ has no $R$-points.

Another natural generalization of Shafarevich's conjecture is to assert that for any global field $K$, the absolute Galois group of the maximal abelian extension of $K$ is free.  In Theorem~4.1 we show that this holds in the function field case.  The proof relies on the corresponding result in [Ha1] and [Po1] about the maximal cyclotomic extension.  The number field case of this conjecture, too, remains open.

While Shafarevich's conjecture and its generalizations above concern one-dimensional function fields, it is also possible to consider higher dimensional versions.  There, for cohomological reasons one must rule out the case of finite coefficient fields; and here we consider function fields over a separably closed field.  In the local case of a smooth surface over a separably closed field $k$, Theorem~4.5 shows that the absolute Galois group of the maximal abelian extension of $k((x,y))$ is free.  The first example of the global case would be to ask whether the same holds for $k(x,y)$ for $k$ separably closed.  This remains open.

\medskip

The key approach in this paper is to use that a profinite group is free if and only if it is projective and quasi-free [HS, Theorem~2.1]; see Section~2 below for definitions.  In this paper it is shown that the absolute Galois group of the function field of any curve over a real closed field is quasi-free, as is the maximal abelian extension of $k((x,y))$ for $k$ any field.  Projectivity is classical in the situation of Theorem~4.2, and it follows from [COP] in that of Theorem~4.5; so in each case freeness then results.
This same approach could also be used to provide another proof of the freeness result of [Ha1] and [Po1] referred to above, and which we use here to obtain the freeness of the maximal abelian extension of the function field of a curve over a finite field.  

One could also use the above approach in considering the two-dimensional global version of Shafarevich's conjecture. Namely, let 
$K\ab$ be the maximal abelian extension of $K = k(x,y)$, for $k$ algebraically closed.  Then the absolute Galois group $G_{K\ab}$ is free if and only if it is projective and quasi-free.  Here $G_{K\ab}$ is the commutator subgroup of $G_K$; and projectivity would follow from knowing for every $\ell$ that a Sylow $\ell$-subgroup of this commutator is a free pro-$\ell$ group.  This condition would imply a conjecture of Bogomolov [Bo], asserting that the commutator subgroup of every Sylow $\ell$-subgroup of $G_K$ is a free pro-$\ell$ group. That conjecture is open, as is the quasi-freeness of $G_{K\ab}$, even in the case of $K = \C(x,y)$.  The above considerations suggest a stronger conjecture that if $k$ is an algebraically closed field (or even a field containing all roots of unity) and $K$ is a function field over $k$, then $G_{K\ab}$ is free.  Such a conjecture has been proposed by F.~Pop.

\medskip

After providing background material and definitions concerning profinite groups, Section~2 of this paper proves (in Theorem~2.4) that the commutator subgroup of a quasi-free group is quasi-free, as well as two related results.  Section~3, which discusses aspects of field arithmetic, shows (in Theorem~3.4) that the absolute Galois group of the function field of a curve over a large field is quasi-free, and that the property of having a free absolute Galois group of countably infinite rank is inherited by abelian extensions (Proposition~3.2).  The main results, on free absolute Galois groups, are shown in Section~4, using results from the previous sections together with other results.

I wish to thank R.~Parimala, J.-L.~Colliot-Th\'el\`ene, M.~Jarden, and F.~Pop for discussions about this material, and MSRI for providing the opportunity to begin working on and discussing these results.

\bigskip

\noindent {\bf Section 2.  Profinite groups.}

\medskip

Let $\Pi$ be a profinite group (i.e.\ an inverse limit of finite groups).  An {\it embedding problem} $\E$ for $\Pi$ is a pair of epimorphisms $\left(\alpha: \Pi \to G, \, f: \Gamma \to G\right)$ of profinite groups; it is {\it non-trivial} if $\ker f$ is non-trivial and it is {\it finite} if $\Gamma$ is finite.  (Here and below, homomorphisms are required to be continuous.)  A {\it weak solution} to $\E = (\alpha, \, f)$ consists of a homomorphism $\lambda:\Pi \to \Gamma$ such that $f \circ \lambda = \alpha$.   A solution is called {\it proper} if it is surjective.  A finite embedding $(\alpha, f)$ in which we also have a splitting $s: G \to \Gamma$ of $f$ is called a {\it finite split embedding problem}.  Every finite split embedding problem has a weak solution given by $s \circ \alpha$.  A profinite group $\Pi$ is {\it projective} if every finite embedding problem for $\Pi$ has a weak solution.  Being projective is equivalent to having cohomological dimension at most $1$ [Gru, Theorem~4] (or [Se, I, \S3.4 Prop.~16 and \S5.9 Prop.~45]).  

A subset $S$ of a profinite group $\Pi$ {\it converges to $1$} if $S \cap (\Pi - N)$ is finite for every open normal subgroup $N$ of $\Pi$. 
Similarly, a map $\phi:S \to G$ to a profinite group $G$ {\it converges to $1$} if $S \cap \phi^{-1}(G - N)$ is finite for every open normal subgroup $N$ of $G$.  The {\it rank} of $\Pi$ is the smallest cardinal number $m = d(\Pi)$ such that $\Pi$ has a set of (topological) generators of cardinality $m$ that converges to $1$.  In fact, if the rank of $\Pi$ is infinite, any two such generating sets have the same cardinality [FJ, Prop.~17.1.2].  Note that if the rank $m$ of $\Pi$ is infinite, then there are at most $m$ (continuous) homomorphisms from $\Pi$ to any finite group, since the kernel must be open.  Thus a finite embedding problem for $\Pi$ can have at most $m$ (weak or proper) solutions, if $m = d(\Pi)$ is infinite.

A profinite group $\Pi$ is {\it free} on a generating set $S$ that converges to $1$ if every map $S \to G$ to a profinite group $G$ that converges to $1$
uniquely extends to a group homomorphism $\Pi \to G$.  For every cardinal $m$ there is a free profinite group of rank $m$ [FJ, \S17.4], denoted $\hat F_m$; this is unique up to isomorphism.  A profinite group $\Pi$ is {\it $\omega$-free} if every finite embedding problem for $\Pi$ has a proper solution.  Every free profinite group is $\omega$-free.  And by a theorem of Iwasawa [Iw, p.567], a profinite group of countable rank is free if and only if it is $\omega$-free.  But this equivalence fails for uncountably generated profinite groups [Ja, Example~3.1].  Instead, there is the following result of Melnikov and Chatzidakis [Ja, Lemma~2.1]: if $m$ is an infinite cardinal, then a profinite group $\Pi$ is free of rank $m$ if and only if every non-trivial finite embedding problem for $\Pi$ has exactly $m$ proper solutions. 

Following [HS] and [RSZ], we say that a profinite group $\Pi$ is {\it quasi-free} if there is a cardinal number $m$ such that every non-trivial finite split embedding problem for $\Pi$ has exactly $m$ proper solutions; to indicate the cardinal, we may say that $\Pi$ is $m$-{\it quasi-free}.  It is easy to see that $m$ is necessarily infinite.  Also, $m$ is necessarily equal to the rank of $\Pi$ [RSZ]; so being $m$-quasi-free is equivalent to being quasi-free of rank $m$.  As a variant on the result of Melnikov and Chatzidakis, if $m$ is an infinite cardinal, then a profinite group $\Pi$ is free of rank $m$ if and only if it is projective and $m$-quasi-free [HS, Theorem~2.1].  If $\Pi$ has countable rank, freeness is also equivalent to the condition that $\Pi$ is projective and every finite split embedding problem for $\Pi$ has a proper solution (in analogy with Iwasawa's theorem) [HS, Corollary~2.8].  

The main goal of this section is to show (in Theorem~2.4 below) that if $\Pi$ is a quasi-free profinite group, then its commutator subgroup $\Pi'$ is also quasi-free, of the same rank.  We begin with some lemmas.

\medskip

\noindent{\bf Lemma 2.1.} {\sl Let $\Pi$ be a profinite group, let $\Pi_1$ be a closed subgroup that contains the commutator subgroup of $\Pi$, and let $\E = ({\alpha:\Pi \to G}, {f:\Gamma \to G})$ be a finite embedding problem for $\Pi$ with proper solutions $\beta_1, \beta_2$ having kernels $M_1, M_2$.  Assume that $Z \cap N = 1$, where $Z$ is the center of $\Gamma$ and $N = \ker\, f$.  Also assume that $\Pi_1 \cap M_1 = \Pi_1 \cap M_2$ and that $\ker\, \alpha \subset \Pi_1M_1$.  Then $M_1 = M_2$.}

\medskip

\noindent{\it Proof.}  Let $\Lambda = \ker\, \alpha$ and let 
$\tilde \Lambda$ be the subgroup of $\Pi$ generated by the kernels $M_1, M_2$ of $\beta_1, \beta_2$.  So $M_i \subset \Lambda$ for $i=1,2$, and hence
$\tilde \Lambda \subset \Lambda$.  Let $N_i = \tilde \Lambda/M_i$ for $i=1,2$.  Thus for $i=1,2$,
we have $|N_i| = (\tilde \Lambda:M_i) = (\Lambda:M_i)/(\Lambda:\tilde \Lambda) = |N|/|\bar N|$  where $\bar N = \Lambda/\tilde \Lambda$.  So $|N_1|=|N_2|$.  
The subgroups $\Pi_1, M_1 \subset \Pi$ are normal (since $\Pi_1$ contains the commutator subgroup of $\Pi$); let $\Lambda_0 \subset \Pi$ be the subgroup they generate.  Here $\Lambda_0$ is generated by $\Pi_1$ and $\Lambda$, since $M_1 \subset \Lambda = \ker\, \alpha \subset \Lambda_0 = \Pi_1 M_1$ by hypothesis.  

Let $\tilde M = \Pi_1 \cap M_1 = \Pi_1 \cap M_2$
and $M_0 =  \Pi_1 \cap \tilde \Lambda$.    Since $\Pi_1 M_1 = \Lambda_0$ and $\Pi_1 M_2 \subset \Pi_1\Lambda \subset \Lambda_0$, we have a natural isomorphism $\Pi_1/\tilde M \iso \Lambda_0/M_1$ and a natural inclusion $\Pi_1/\tilde M \hookrightarrow \Lambda_0/M_2$.  These respectively restrict to an isomorphism $N_0 := M_0/\tilde M \iso \tilde \Lambda/M_1 =  N_1$ and to an inclusion $N_0 = M_0/\tilde M \hookrightarrow \tilde \Lambda/M_2 = N_2$.  But $|N_1| = |N_2|$; so this last inclusion must also be an isomorphism.  Hence $N_0, N_1, N_2$ are isomorphic.  Also, for $i=1,2$, the isomorphism $M_0/\tilde M \iso \tilde \Lambda/M_i$ shows that $M_0,M_i$ generate $\tilde \Lambda$.

Let $\Pi_0 = \langle \Pi_1, M_1 \cap M_2 \rangle$.  So the natural map $\Pi_0/(M_1 \cap M_2) \to \Pi_1/(M_1 \cap M_2 \cap \Pi_1) = \Pi_1/\tilde M$ is an isomorphism.  Also, $\Pi_0 \cap M_i$ and $\Pi_1$ generate $\Pi_0$ for $i=1,2$, since $M_1 \cap M_2 \subset \Pi_0 \cap  M_i\subset \Pi_0$.  So the natural inclusion 
$\Pi_0/(\Pi_0 \cap M_i) \to \Pi_1/(\Pi_1 \cap M_i) = \Pi_1/\tilde M$ is also an isomorphism, for $i=1,2$.  But $M_1 \cap M_2 \subset \Pi_0 \cap M_i$;
so $\Pi_0 \cap M_1 = M_1 \cap M_2 = \Pi_0 \cap M_2$.  Moreover  $\Pi_0M_1 = \langle \Pi_1, M_1 \cap M_2, M_1 \rangle = \Pi_1 M_1$ contains $\Lambda = \ker\, \alpha$, by hypothesis.  Also, $\Pi_0$ contains the commutator subgroup of $\Pi$, since it contains $\Pi_1$.  So it suffices to prove the lemma with $\Pi_1$
replaced by $\Pi_0$.  Thus we may assume that $\Pi_1$ contains $M_1 \cap M_2$, and hence that $M_1 \cap M_2 = \Pi_1 \cap M_1 \cap M_2 = \Pi_1 \cap M_1 = \Pi_1 \cap M_2 = \tilde M$.  

So $M_1, M_2$ are normal subgroups of $\tilde \Lambda$ satisfying $M_1 \cap M_2 = \tilde M$ and $M_1M_2 = \tilde \Lambda$; thus $\tilde \Lambda/\tilde M \approx \tilde \Lambda/M_1 \times \tilde \Lambda/M_2 \approx N_1 \times N_2$, with $M_1/\tilde M \approx 1 \times N_2$ and $M_2/\tilde M \approx N_1 \times 1$ under this isomorphism.  The subgroup $M_0 = \Pi_1 \cap \tilde \Lambda\subset\tilde \Lambda$ is normal; let $Q$ be the quotient.  Thus we have an exact sequence $1 \to N_0 \to N_1 \times N_2 \to Q \to 1$.  Now $\tilde M = \Pi_1 \cap M_2 = M_0 \cap M_2$ (using $\tilde M \subset M_0 \subset \Pi_1$); so $N_0 \cap (1 \times N_2) = 1$.  Also, $N_0$ and $1 \times N_2$ are normal subgroups of $N_1 \times N_2$ that generate $N_1 \times N_2$, because $M_0$ and $M_1$ generate $\tilde \Lambda$.
So $N_1 \times N_2$ is isomorphic to $N_0 \times (1 \times N_2)$, and $Q = (N_1 \times N_2)/N_0 \approx 1 \times N_2 \approx N_2 \approx N_1$.  Since $M_0 = \Pi_1 \cap \tilde \Lambda$, we have that $Q = \tilde \Lambda/M_0$ is isomorphic to a subgroup of $\Pi/\Pi_1$, and so is abelian.  Hence so are the isomorphic groups $N_i$, for $i=0,1,2$.  

Since $\Pi_1$ contains the commutator subgroup of $\Pi$, the group $\Pi_1/\tilde M$ contains the commutator subgroup of $\tilde \Gamma := \Pi/\tilde M$ and hence contains $[\tilde \Gamma,M_2/\tilde M]$.  But since $M_2/\tilde M$ is normal in $\tilde \Gamma$, we also have $[\tilde \Gamma,M_2/\tilde M] \subset M_2/\tilde M$.  So in fact we have the containment $[\tilde \Gamma,M_2/\tilde M] \subset (\Pi_1/\tilde M) \cap (M_2/\tilde M) = 1$ since $\Pi_1\cap M_2 = \tilde M$.  That is, $M_2/\tilde M = N_1 \times 1$ is in the center of $\tilde \Gamma$.  Taking images under $1 \times N_2$, i.e.\ under $\tilde \Gamma \to \Gamma = \Pi/M_1$, we have that the subgroup $N_1 = \tilde \Lambda/M_1$ of $N = \Lambda/M_1$
is in the center $Z$ of $\Gamma = \Pi/M_1$.  But $Z \cap N = 1$ by hypothesis; so $N_1=1$ and hence also $N_2=1$. Thus $M_1 = \tilde\Lambda = M_2$, as required.
\qed

\medskip

\noindent{\bf Lemma~2.2.} {\sl Let $\Pi$ be a quasi-free profinite group of rank $m$, and let $\Pi'$ be its commutator subgroup.  Let $p$ be a prime number.  Then there are $m$ closed normal subgroups of $\Pi'$ having index $p$.}

\medskip

\noindent {\it Proof.}  Define a $2 \times 2$-matrix $A_p$ over $\F_p$ as follows: If $p=2$ the rows are $(0 \ 1)$ and $(1 \ 1)$; and if $p \ne 2$ the rows are $(0 \ a)$ and $(1 \ 0)$, where $a \in \F_p^*$ is not a square.  Let $r$ be the order of $A_p$ in ${\rm GL}(2,p)$.  So the cyclic group $C_r$ acts on the two-dimensional $\F_p$-vector space $\F_p^2$ via left multiplication of the matrix $A_p$ on column vectors.  This action is irreducible over $\F_p$, since the minimal polynomial $f_p(x) \in \F_p[x]$ of $A_p$ is irreducible.  (Namely, $f_2(x) = x^2-x-1$ and $f_p(x) = x^2-a$ for $p$ odd.) 

Let $\Gamma = C_p^2 \sdp C_r$, where the conjugation action of $C_r$ on the group $C_p^2 \approx \F_p^2$ is via $A_p$ as above.  Since this action is irreducible, $C_p^2$ has no non-trivial proper subgroup that is normal in $\Gamma$.  Note that if $Z$ is the center of $\Gamma$, then $C_p^2$ is not contained in $Z$; and hence $Z \cap C_p^2 = 1$, being a normal subgroup of $\Gamma$ contained in $C_p^2$.  Also, 
the commutator subgroup $\Gamma'$ of $\Gamma$ is non-trivial since $\Gamma$ is non-abelian, and it is contained in $C_p^2$ since $\Gamma/C_p^2$ is abelian.  So again by irreducibility, $\Gamma' = C_p^2$.  That is, the maximal abelian quotient of $\Gamma$ is $C_r$.  

Since $\Pi$ is quasi-free of rank $m$, $\Pi$ has $m$ distinct normal subgroups with quotient group isomorphic to $C_r$.  Picking one of them, say $\Lambda$, and an epimorphism $\alpha:\Pi \to C_r$ with kernel $\Lambda$, we obtain a finite split embedding problem for $\Pi$ given by this map and the split exact sequence $1 \to C_p^2 \to \Gamma \to C_r \to 1$.  Again, since $\Pi$ is quasi-free of rank $m$, there are $m$ distinct proper solutions $\beta:\Pi \to \Gamma$ to this embedding problem.  Taking kernels, we have that there are $m$ distinct normal subgroups $M \subset \Pi$ such that $\Pi/M \approx \Gamma$ and $M \subset \Lambda$ (using that there are only finitely many epimorphisms to a given finite group with a given kernel).  For each such $M$, the subgroup of $\Pi$ generated by $\Pi'$ and $M$ is the minimal normal subgroup of $\Pi$ that contains $M$ and has abelian quotient;
i.e.\ $\Pi'M = \Lambda = \ker\, \alpha$ (by the maximality assertion in the previous paragraph).  Thus $\Pi'/(\Pi' \cap M) \approx \Lambda/M = \ker((\Pi/M)\to(\Pi/\Lambda)) \approx \ker(\Gamma \to C_r)
= C_p^2$.  

So for any two such normal subgroups $M_1, M_2$ (among the $m$ given by solutions to the embedding problem), we may apply Lemma~2.1; and conclude that if $\Pi' \cap M_1 = \Pi' \cap M_2$ then $M_1=M_2$.  That is, the $m$ solutions to the embedding problem induce $m$ distinct normal subgroups of $\Pi'$ having quotient $C_p^2$.  Taking the inverse images of $1 \times C_p$ and of $C_p \times 1$ under such quotient maps, we obtain $m$ normal subgroups of $\Pi'$ having quotient group $C_p$, possibly with repetitions.  If there are precisely $m' \le m$ distinct normal subgroups of $\Pi'$ with quotient $C_p$ arising this way, then the number of normal subgroups of $\Pi'$ with quotient $C_p^2$ obtained by taking intersections is also $m'$.  But there are (at least) $m$ of them, as noted above.  So in fact $m'=m$.
\qed

\medskip

\noindent {\it Remark.} Lemma 2.2 is a weak form of Theorem~2.4 below, in the case of a split embedding problem corresponding to a short exact sequence $1 \to C_p \to C_p \to 1 \to 1$.  Note that Lemma~2.2 did not assert that there are {\it no more} than $m$ solutions; but this is in fact true and follows from Theorem~2.4.

\medskip

Let $\iota:\Pi_0 \to \Pi_1$ be a homomorphism of profinite groups (e.g.\ an inclusion), let $f:\Gamma \to G$ be an epimorphism of finite groups, and let $\alpha_i:\Pi_i \to G$ be an epimorphism for $i=0,1$.  Thus 
$\E_i = (\alpha_i:\Pi_i \to G, f:\Gamma \to G)$ is a finite embedding problem for $\Pi_i$.  We say that $\E_1$ {\it induces} $\E_0$ if $\alpha_0 =  \alpha_1 \circ \iota$.  If $\E_1$ induces $\E_0$ and if
$\beta_i:\Pi_i \to \Gamma$ is a weak solution to $\E_i$ for $i=0,1$, we say that $\beta_1$ {\it induces} $\beta_0$ if $\beta_0 = \beta_1 \circ \iota$.  Note that if $\beta_1$ is a proper solution to $\E_1$, then the induced solution $\beta_0$ to $\E_0$ need not be proper.

\medskip

\noindent{\bf Lemma 2.3.} {\sl Let $\Pi$ be a profinite group, let $\Pi'$ be its commutator subgroup, and let $\E = (\alpha:\Pi' \to G, f:\Gamma \to G)$ be a non-trivial finite split embedding problem for $\Pi'$.

a) Then there is a finite index closed normal subgroup $\Pi_1 \
\subset\Pi$ containing $\Pi'$ together with an embedding problem $\E_1 =  ({\alpha_1:\Pi_1 \to G}, {f:\Gamma \to G})$ that induces $\E$. 

b) If $\Pi$ is quasi-free of rank $m$ then $\Pi_1$ and $\E_1$ can be chosen such that $\E_1$ has a set of $m$ proper solutions each of which induces a proper solution to $\E$.}

\medskip

\noindent{\it Proof.}  
Let $\tilde \Lambda = \ker\, \alpha \subset \Pi'$ and let $N = \ker\, f \ne 1$.

a) Since $\alpha$ is continuous, $\tilde \Lambda$ is a closed normal subgroup of $\Pi'$, where $\Pi'$ has the topology induced from that of the profinite group $\Pi$.  Hence there is a finite index closed subgroup $\Lambda_1 \subset \Pi$ such that $\Lambda_1 \cap \Pi' = \tilde \Lambda$.  Let $\Pi_1 \subset \Pi$ be the subgroup generated by $\Lambda_1$ and $\Pi'$.  So $\Pi'$ is a closed normal subgroup of $\Pi_1$, since it is a closed normal subgroup of $\Pi$.  Since $\Pi'$ and $\Lambda_1$ generate $\Pi_1$ and have intersection $\tilde \Lambda$, and since $\tilde \Lambda$ is normal in $\Pi'$,
it follows that $\Lambda_1$ is normal in $\Pi_1$ and
the natural inclusion $G = \Pi'/\tilde \Lambda \hookrightarrow \Pi_1/\Lambda_1$ is an isomorphism.  So we obtain an epimorphism $\alpha_1:\Pi_1 \to G$ whose restriction to $\Pi'$ is $\alpha$.  Thus the embedding problem $\E_1 =  (\alpha_1:\Pi_1 \to G, f:\Gamma \to G)$ induces $\E$, as desired.

b) Let $\Pi_0 \subset \Pi$ be a finite index closed normal subgroup containing $\Pi'$, together with an embedding problem $\E_0 =  (\alpha_0:\Pi_0 \to G, f:\Gamma \to G)$ that induces $\E$; these exist by part (a).  Let $\Lambda_0 = \ker\, \alpha_0 \subset \Pi_0$.  So $\tilde \Lambda = \Pi' \cap \Lambda_0$.

Since $\Pi$ is quasi-free of rank $m$, there are $m$ closed normal subgroups $\Phi \subset \Pi$ such that $\Pi/\Phi \approx C_2$ (arising from the embedding problem for $\Pi$ corresponding to the exact sequence $1 \to C_2 \to C_2 \to 1 \to 1$).  Thus in particular we may choose such a $\Phi$ that does not contain $\Lambda_0$ (since only finitely many subgroups of $\Pi$ contain the finite index subgroup $\Lambda_0$).  Note that $\Phi$ and $\Lambda_0$ generate $\Pi$, since $(\Pi:\Phi)=2$.  Let $\Pi_1 = \Pi_0  \cap \Phi \subset \Pi$ and let $\Lambda_1 = \Lambda_0\cap\Phi = \Lambda_0 \cap \Pi_1
\subset \Pi$.  Thus $\Pi_1$ is a closed normal subgroup of $\Pi_0$, and $\Lambda_1$ is a closed normal subgroup of  the groups $\Pi_1$, $\Lambda_0$ and $\Pi_0$.  Here $(\Lambda_0:\Lambda_1) = (\Pi_0:\Pi_1) = 2$, since $(\Pi:\Phi)=2$ and since $\Phi$ does not contain $\Lambda_0$.  So $\Pi_0$ is generated by $\Lambda_0$ and $\Pi_1$.  Hence the natural map $\Pi_1/\Lambda_1 \hookrightarrow \Pi_0/\Lambda_0 = G$ (through which the isomorphism $G = \Pi'/\tilde\Lambda \iso \Pi_0/\Lambda_0$ factors) is an isomorphism, and we have isomorphisms $\Pi_0/\Lambda_1 \iso \Pi_0/\Lambda_0 \times \Pi_0/\Pi_1 \iso G \times C_2$.  So the restriction $\alpha_1:\Pi_1 \to G$ of $\alpha_0:\Pi_0 \to G$ is surjective with kernel $\Lambda_1$, and $\alpha_1$ in turn restricts to $\alpha:\Pi' \to G$, whose kernel $\tilde\Lambda$ is contained in $\Lambda_1$.  Thus the embedding problem $\E_1 =  (\alpha_1:\Pi_1 \to G, f:\Gamma \to G)$ induces $\E$.  
So the natural map $G = \Pi'/\tilde \Lambda \to \Pi_1/\Lambda_1$ is an isomorphism; and hence $\Pi_1 = \Pi' \Lambda_1$.
Moreover $\alpha_1$ lifts to a surjection $\hat\alpha_1:\Pi_0 \to G \times C_2$ having kernel $\Lambda_1$, corresponding to the above isomorphism $\Pi_0/\Lambda_1 \iso G \times C_2$.   

Since the closed subgroup $\Pi_0 \subset \Pi$ has finite index, it is also open.  But every open subgroup of a quasi-free group is also quasi-free of the same rank [RSZ].  So $\Pi_0$ is quasi-free of rank $m$.  
Let $\hat\Gamma$ be the semi-direct product of $N \times N$ with the group $G \times C_2$, where $G$ acts on each factor $N$ as it does in $\Gamma$, and where $C_2$ acts by interchanging the two copies of $N$. Also, let $\hat f_1: \hat \Gamma \to G \times C_2$ be the canonical surjection, and consider the finite split embedding problem $\hat\E_1 = (\hat\alpha_1:\Pi_0 \to G \times C_2,\ \hat f_1:\hat\Gamma \to G \times C_2)$ for $\Pi_0$.  Since $\Pi_0$ is quasi-free of rank $m$, this embedding problem has $m$ proper solutions.

Consider any proper solution to $\hat \E_1$, say $\hat \beta_1:\Pi_0 \to \hat\Gamma$.  So $M := \ker\, \hat \beta_1$ is normal in $\Lambda_1$ and in $\Pi_0$, with quotient groups $\Lambda_1/M \approx N \times N$ and
$\Pi_0/M \approx \hat \Gamma$.  Here $H := \Lambda_0/M = (N \times N) \sdp C_2$, with $C_2 = \Lambda_0/\Lambda_1$ interchanging the two factors of $N \times N = \Lambda_1/M$.  Let $M_1$ be the inverse image of $1 \times N$ under the quotient map $\Lambda_1 \to N \times N$; thus $\Pi_1/M_1 \approx \Gamma$.  Also, $M$ is the largest normal subgroup of $\Pi_0$ that is contained in $M_1$ (since any such subgroup would also have to be contained in the inverse image of $N \times 1$); so $M$ is determined by $M_1$ and thus distinct choices of $M$ lead to distinct choices of $M_1$.  Thus there are $m$ distinct choices for $M_1$, arising from the $m$ choices for $M$.  Each such choice for $M_1$ is the kernel of a proper solution $\beta_1$ to the embedding problem $\E_1$, inducing a weak solution $\beta := \beta_1|_{\Pi'}:\Pi' \to \Gamma$ to $\E$ with kernel $\Pi' \cap M_1$.  It suffices to show that in fact $\beta$ is a proper solution to $\E$, i.e.\ is surjective.

If $n\in N$ and $\iota$ is the involution in $C_2$, the commutator $[(n,1),\iota] \in [N \times N, C_2] \subset H = (N \times N) \sdp C_2$ is equal to $(n,n^{-1}) \in N \times N$.  
Thus $N \times N = \Lambda_1/M$ is generated by $1 \times N = M_1/M$ and the commutator subgroup $H'$ of $H = \Lambda_0/M$ (where $H' \subset N \times N$ because $H/(N \times N)$ is abelian).  So $\Lambda_1$ is generated by $M_1$ and $\Lambda_0'$, the commutator subgroup of $\Lambda_0$.    Since $\Lambda_0' \subset \Pi'$, we have that 
$\Pi'M_1 = \langle \Pi', \Lambda_0', M_1 \rangle = \Pi' \Lambda_1 = \Pi_1$.  Hence the 
natural inclusion $\beta(\Pi') \approx \Pi'/\ker\, \beta = \Pi'/(\Pi' \cap M_1) \hookrightarrow \Pi_1/M_1= \Gamma$ is an isomorphism.   So $\beta:\Pi' \to \Gamma$ is surjective, as desired. \qed

\medskip

\noindent{\bf Theorem 2.4.} {\sl Let $m$ be an infinite cardinal.  If $\Pi$ is a quasi-free profinite group of rank $m$, then so is its commutator subgroup $\Pi'$.}

\medskip

\noindent{\it Proof.}  Let $\E = (\alpha:\Pi' \to G, f:\Gamma \to G)$ be any non-trivial finite split embedding problem  for $\Pi'$.  Let $N = \ker\, f \ne 1$ and let $\tilde \Lambda = \ker\, \alpha$.  Thus $\Pi'/\tilde \Lambda \approx G$.
We wish to show that $\E$ has exactly $m$ distinct proper solutions.

Suppose that $\E$ has more than $m$ proper solutions $\beta:\Pi' \to \Gamma$.  For each $\beta$, there is a finite index closed normal subgroup $\Pi_0 \subset \Pi$ containing $\Pi'$ and an epimorphism $\beta_0:\Pi_0 \to \Gamma$ such that $\beta = \beta_0 \circ \iota$, where $\iota:\Pi' \hookrightarrow \Pi_0$ is the containment map.  (Namely, this follows from Lemma~2.3(a), applied to the trivial embedding problem $(\beta:\Pi' \to \Gamma, {\rm id}:\Gamma \to \Gamma)$.)  
Now $\Pi_0$ has finite index and contains $\Pi'$; so $\Pi_0$ is the kernel of an epimorphism $\kappa:\Pi \to D$ for some finite abelian group $D$.  But
there are only countably many finite abelian groups $D$ up to isomorphism; and for each $D$ there are exactly $m$ epimorphisms $\Pi \to D$ (since $\Pi$ is quasi-free of rank $m$).  So by the assumption that $\E$ has more than $m$ proper solutions, there must be some finite abelian group $D$ and some epimorphism $\kappa_0:\Pi \to D$, say with kernel $\Pi_0$, such that there are more than $m$ epimorphisms $\beta_0:\Pi_0 \to \Gamma$.  But $\Pi_0$ is an open subgroup of $\Pi$ (being closed and of finite index), and hence is quasi-free of rank $m$ by [RSZ]. This is a contradiction.  So actually $\E$ has at most $m$ proper solutions. 

It remains to show that $\E$ has at least (and hence exactly) $m$ proper solutions.  

By Lemma~2.3, there is a finite index closed normal subgroup $\bar\Pi$ of $\Pi$ containing $\Pi'$, together with an embedding problem 
$\bar\E = ({\bar\alpha:\bar\Pi \to G}, {f:\Gamma \to G})$ that induces $\E$, such that $\bar\E$ has $m$ proper solutions each of which induces a proper solution to $\E$.  Say that $\alpha_1$ is such a proper solution of $\bar\E$, with kernel $M_1 \subset \bar\Pi$.
Since $\bar\E$ induces $\E$, we have that $\tilde\Lambda = \ker\, \alpha = \Pi' \cap M_1$ and the map $\Pi'/(\Pi' \cap M_1) \to \bar\Pi/M_1$ is an isomorphism.  So $\Pi'M_1=\bar\Pi \supset \ker\, \alpha_1$.  Also, $\Pi'$ contains the commutator subgroup of $\bar\Pi$.  Let $Z$ be the center of $\Gamma$.  By Lemma~2.1 (with $\bar \Pi$ and $\Pi'$ here playing the roles of $\Pi$ and $\Pi_1$ there), if $Z \cap N = 1$ then distinct choices of $\mu_1$ that have distinct kernels $M_1$ yield distinct proper solutions to $\E$.  So in this case we are done; and we are therefore reduced to the case that $Z \cap N \ne 1$. 

We may thus assume that there is a cyclic subgroup $C$ of prime order $p$ in ${Z \cap N}$. According to Lemma~2.2, there are (at least) $m$ distinct closed normal subgroups of index $p$; and so $\Hom(\Pi',C)$ has cardinality $\ge m$.  Let $\beta:\Pi' \to \Gamma$ be the proper solution to $\E$ given by {\it some} choice of $\alpha_1$ and $M_1$ in the previous paragraph.
Since $C$ is central in $\Gamma$, for each $\varepsilon \in \Hom(\Pi',C)$ we obtain a homomorphism $\beta\cdot \varepsilon:\Pi' \to \Gamma$ given by $(\beta\cdot \varepsilon)(a) = \beta(a)\varepsilon(a)$ for $a \in \Pi'$.  The composition of $\beta\cdot \varepsilon$ with the quotient map $\Gamma \to G$ is the surjection $\alpha:\Pi' \to G$, since this is true for $\beta$ and since $C \subset N = \ker(\Gamma \to G)$.  Moreover, the compositions of $\beta$ and of $\beta \cdot \varepsilon$ with $\Gamma \to \Gamma/C$ also agree, and the former is surjective; so $\Gamma$ is generated by $C$ and the image of $\beta\cdot \varepsilon$.  In particular, $\beta\cdot \varepsilon$ is surjective if and only if its image contains $C$.  Also, distinct choices of $\varepsilon$ yield distinct homomorphisms $\beta\cdot \varepsilon$.  So it suffices to show that the image of $\beta\cdot \varepsilon$ contains $C$ for at least $m$ choices of $\varepsilon \in \Hom(\Pi',C)$.

Let $\Delta \subset \Pi'$ be the inverse image of $C$ under $\beta$.
Since $\beta:\Pi'\to \Gamma$ is surjective, the image of $\beta|_\Delta$ is $C$.  Also, for each $\varepsilon \in \Hom(\Pi',C)$, the image of the restriction $(\beta\cdot \varepsilon)|_\Delta$ is either $C$ or $1$.
Let $S$ be the set of $\varepsilon \in \Hom(\Pi',C)$ such that this image is $C$.  For any $\varepsilon \in S$, the map $\beta\cdot \varepsilon$ is surjective, since its image contains $C$; and so it suffices to show that the cardinality of $S$ is at least $m$.   If the complement of $S$ in $\Hom(\Pi',C)$ has cardinality less than $m$, then the cardinality of $S$ is at least $m$, and we are done.  On the other hand, if the cardinality of the complement of $S$ is $\ge m$, then fix some $\varepsilon_0$ in this complement.  For any {\it other} $\varepsilon$ in the complement of $S$, consider the map $u_\varepsilon := \varepsilon_0^{-1}\cdot \varepsilon:\Pi' \to C$ sending $a \in \Pi'$ to $\varepsilon_0(a)^{-1}\varepsilon(a)$.  The restriction of this map to $\Delta$ is trivial, since $(\beta\cdot \varepsilon)(a) = (\beta\cdot \varepsilon_0)(a) = 1$ for $a \in \Delta$. So $(\beta\cdot u_\varepsilon)|_\Delta = \beta|_\Delta$, whose image is $C$; and hence $u_\varepsilon \in S$.  Since distinct $\varepsilon$'s in the complement of $S$ induce distinct $u_\varepsilon$'s in $S$, it follows that the cardinality of $S$ is at least $m$.
\qed

\medskip

\noindent{\bf Corollary 2.5.} {\sl Let $m$ be an infinite cardinal.  If $\Pi$ is a free profinite group of rank $m$, then so is its commutator subgroup $\Pi'$.}

\medskip

\noindent{\it Proof.} Since $\Pi$ is free profinite of rank $m$, and since $m$ is infinite, we have that $\Pi$ is quasi-free of rank $m$.  So by Theorem~2.4, $\Pi'$ is also quasi-free of rank $m$.  Since $\Pi'$ is a closed subgroup of the free profinite group $\Pi$, it follows from [FJ, Corollary~22.4.6] that $\Pi'$ is projective.  Since $\Pi'$ is projective and quasi-free of rank $m$, it is free of rank $m$ by [HS, Theorem~2.1]. \qed

\medskip

In the case that $m$ is countable, a stronger conclusion is possible (see Proposition~3.2(b) below).

The following result is another variant of Theorem~2.4, considering just the existence of finite quotients rather than embedding problems.  

\medskip

\noindent{\bf Proposition 2.6.} {\sl Let $\Pi$ be a profinite group with the property that every finite group is a quotient of $\Pi$ by a closed normal subgroup of finite index.  Then the commutator subgroup $\Pi'$ of $\Pi$ also has this property.}

\medskip

\noindent{\it Proof.}  We proceed as at the end of the proof of Lemma~2.3.  Let $N$ be any finite group, and let $H = (N \times N) \sdp C_2$, where $C_2$ acts by interchanging the two copies of $N$.  By hypothesis, $\Pi$ has a closed normal subgroup $M$ such that $\Pi/M = H$.  Let $p:\Pi \to H$ be the canonical surjection, let $\Pi_1 = 
p^{-1}(N \times N)$, and let $M_1 = p^{-1}(1 \times N)$.  Thus $M \subset M_1 \subset \Pi_1 \subset \Pi$, and $\Pi' \subset \Pi_1$ since $\Pi/\Pi_1$ is abelian.  As in the proof of Lemma~2.3, $M$ is the largest normal subgroup of $\Pi$ contained in $M_1$, and $N \times N = \Pi_1/M$ is generated by $1 \times N = M_1/M$ and the commutator subgroup $H'$ of $H = \Pi/M$.  So $\Pi_1$ is generated by $M_1$ and $\Pi'$.  Hence the natural inclusion $\Pi'/(\Pi' \cap M_1) \hookrightarrow  \Pi_1/M_1 \approx N$ is an isomorphism.  Thus $\Pi' \cap M_1$ is a closed normal subgroup of $\Pi'$ with quotient group isomorphic to $N$. \qed

\bigskip

\noindent {\bf Section 3. Field arithmetic.}

\medskip

Let $K$ be a field, with separable closure $K\s$.  The {\it absolute Galois group} of $K$ is the profinite group $G_K := \Gal(K\s/K)$.  An {\it embedding problem} for $K$ is an embedding problem $\E = \left(\alpha: G_K \to G, \, f: \Gamma \to G\right)$ for $G_K$.  Here the epimorphism $\alpha$ corresponds to a $G$-Galois field extension $L$ of $K$ together with a $K$-inclusion $i:L \hookrightarrow K\s$.  
A proper solution to $\E$ corresponds to a $\Gamma$-Galois field extension $M$ of $K$ that contains $L$, together with a $K$-inclusion $j:M \hookrightarrow K\s$ that extends $i$.  Thus $\E$ has a proper solution if and only if the given $G$-Galois field extension of $K$ can be embedded into a $\Gamma$-Galois field extension (hence the terminology).  Note that if $G$ and $\Gamma$ are finite, then there are only finitely many $K$-inclusions $i$ and $j$ as above, for given field extensions $L$ and $M$.  Also, if $m := \card K$ is infinite, then $K$ has at most $m$ field extensions of finite degree; and so a finite embedding problem for $K$ can have at most $m$ (weak or proper) solutions.

In the above situation, suppose that $K$ is a function field over a subfield $F$ (i.e.\ of finite transcendence degree over $F$, with $F$ algebraically closed in $K$), and let $\beta$ be a proper solution to $\E = \left(\alpha: G_K \to G, \, f: \Gamma \to G\right)$ corresponding to a pair $(M,j)$ extending $(L,i)$.  We say that the proper solution $\beta$ is {\it regular} (with respect to $F$) if the algebraic closures of $F$ in $L$ and in $M$ are the same (regarding $L \subset M$).

The Galois cohomology of a field $K$ is the same as the group cohomology of $G_K$, and so $K$ and $G_K$ have the same cohomological dimension.  We say that $K$ is {\it free} [resp.\ {\it quasi-free}, {\it $\omega$-free}, {\it projective}] if $G_K$ is.  So $K$ is projective if and only if it has cohomological dimension $\le 1$.  Also, if $K$ is quasi-free of rank $m_0$, then $\card K \ge m_0$.  
We say that a profinite group $G$ {\it is a Galois group over} $K$ if there is a Galois field extension $L$ of $K$ with Galois group isomorphic to $G$; this is equivalent to saying that $G_K$ has a closed normal subgroup $N$ such that $G_K/N$ is isomorphic to $G$.

For any field $K$, let $K\ab$ denote its maximal abelian extension (in a given separable closure).  By considering the absolute Galois group $\Pi = G_K$ and its commutator $\Pi' = G_{K\ab}$, we may restate 
Proposition~2.6, Theorem~2.4 and Corollary~2.5 in field-theoretic terms as follows:

\medskip

\noindent{\bf Proposition 3.1.} {\sl Let $K$ be a field.

a) If $K$ has the inverse Galois property (i.e.\ every finite group is a Galois group over $K$), then the same holds for $K\ab$.

b) Let $m$ be an infinite cardinal.  If $K$ is quasi-free [resp.\ free] of rank $m$, then so is $K\ab$.}

\medskip

Recall that a field $K$ is {\it large} [Po2] (also sometimes called {\it ample} [FJ, Remark~16.12.3]) if every smooth $K$-curve with a $K$-point has infinitely many $K$-points.  Examples of large fields include fraction fields of henselian (e.g.\ complete) discrete valuation rings; real closed fields (e.g.\ $\R$); the field of totally real (or totally $p$-adic) algebraic numbers; algebraically closed fields; more generally pseudo-algebraically closed fields (PAC fields: fields $K$ such that smooth geometrically integral $K$-variety has a $K$-rational point); and algebraic extensions of large fields [Po2].  The property of being large is equivalent to the property that for every smooth integral $K$-variety $X$, if $X$ has a $K$-point then $X(K)$ is Zariski dense (using that the union of smooth $K$-curves containing a given smooth $K$-point on an integral $K$-variety $X$ is Zariski dense in $X$).  It is also equivalent to the condition that $K$ is existentially closed in $K((t))$; i.e.\ every $K$-variety with a $K((t))$ point has a $K$-point.  See [Po2, Proposition~1.1].

A key property of large fields (first shown by F.~Pop) is the following:  Let $F$ be a large field and let $K$ be the function field of a smooth projective $F$-curve.  Then every finite split embedding problem for $K$ has a proper regular solution.  Versions of this result have appeared in [Po1], [Po2], and [HJ] (and see [Ha3, \S5.1] for a further discussion).  Hence large Hilbertian fields $K$ 
have the property that every finite split embedding problem has a proper solution.  If in addition $G_K$ is projective then $G_K$ is $\omega$-free; and if also $K$ is countable then Iwasawa's theorem [Iw, p.567] applies and so $G_K$ is free of countable rank [Po2, Theorem~2.1].

Using these ideas, we obtain a stronger form of the free case of Proposition~3.1(b) and of Corollary~2.5 when $m$ is countable.  We state this both in the group-theoretic and field-theoretic settings.  

\medskip

\noindent{\bf Proposition 3.2.} {\sl 
a) Let $\Pi$ be a profinite group and let $\Pi_1$ be a closed subgroup of $\Pi$ that contains the commutator subgroup $\Pi'$ of $\Pi$.  If $\Pi$ is free of countably infinite rank, then so is $\Pi_1$. 

b) Let $K_1$ be an abelian extension of a field $K$.  If the absolute Galois group of $K$ is free of countably infinite rank, then the same holds for $K_1$.}

\medskip

\noindent{\it Proof.}  First consider (b) in the special case that $K$ is a countable PAC field.  Since the absolute Galois group $G_K$ is free of infinite rank, it is $\omega$-free.  Being PAC and $\omega$-free implies that $K$ is Hilbertian, by a theorem of Roquette [FJ, Corollary~27.3.3].  Any algebraic extension of a PAC field is PAC [FJ, Corollary~11.2.5], and any abelian extension of a Hilbertian field is Hilbertian [FJ, Theorem~16.11.3]. So $K_1$ is also a countable Hilbertian PAC field.  
Since $K_1$ is PAC, it is large and its absolute Galois group $G_{K_1}$ is projective [FJ, Corollary~23.1.3] (with projectivity also following from the fact that $G_{K_1}$ is a closed subgroup of the free profinite group $G_K$ [FJ, Corollary~22.4.6]).  So $K$ is countable, large and Hilbertian, with $G_K$ projective; hence [Po2, Theorem~2.1] applies and asserts that $G_K$ is free of countable rank (see also [FJ, Example~24.8.5].)  This shows that (b) holds in this special case.

Next, observe that (a) follows from the above special case of (b) by using the fact that a free profinite group of countably infinite rank is isomorphic to the absolute Galois group of any countable Hilbertian PAC field $K$, and by taking $K_1$ to be the fixed field of $\Pi_1$ in the separable closure of $K$. 

Finally, observe that the general case of (b) follows from (a) by letting $\Pi$ and $\Pi_1$ be the absolute Galois groups of $K$ and $K_1$ respectively. \qed 

\medskip

As in [HS], call a field $K$ {\it very large} if every smooth $K$-curve with a $K$-point has exactly $m$ $K$-points, where $m$ is the cardinality of $K$.  This is equivalent to the property that for every smooth integral $K$-variety $X$, if $X$ has a $K$-point then every non-empty open subset of $X$ contains exactly $m$ $K$-points (using the same reasoning as for the corresponding characterization of large).

Observe that every large field is infinite, as is every very large field (e.g.\ by considering the curve $\P^1_K$).  Hence every very large field is large.  Also, if $K$ is an infinite field of cardinality $m$, then every $K$-variety (of finite type) has at most $m$ $K$-points. 

The proof of the following proposition is due to F.~Pop.

\medskip

\noindent{\bf Proposition 3.3.} [Pop] {\sl Let $K$ be a large field of cardinality $m$.  Then $K$ is very large.} 

\medskip

\noindent{\it Proof.}  Let $X$ be a smooth $K$-curve with a $K$-point $P$, where $K$ is large.  We wish to show that the cardinality of $X(K)$ is equal to $m$.  Since $X$ is a $K$-variety, $X(K)$ has cardinality at most $m=\card K$.  So it suffices to prove the reverse inequality; and for this we may assume that $X$ is connected.  Possibly after deleting finitely many points (other than $P$) from $X$, we may embed $X$ in $\A^2_K$.  After replacing $X$ by its image in $\A^2_k$, and making a change of variables in the plane, we may assume that $X$ is a smooth plane curve containing the origin, defined by a polynomial $f$ such that $\partial f/\partial y$ does not vanish at the origin.  We claim that for each $a \in K$ we may choose a pair of $K$-points $(x_1,y_1),(x_2,y_2) \in X(K)$ such that $x_2 \ne 0$ and $x_1/x_2 = a$.  If this is shown, we obtain an injection $i:K \hookrightarrow X(K) \times X(K)$; and this then implies that the cardinality of $X(K)$ is at least $m$, as desired.

So it suffices to prove the claim.  Let $a \in K$.  Consider affine 4-space $\A^4_K$ with coordinates $X_1,Y_1,X_2,Y_2$, and the subvariety $V_a\subset \A^4_k$ defined by:
$$f(X_1,Y_1)=0, \ f(X_2,Y_2)=0, \ X_1-aX_2=0.$$
Here $V_a\cong(X\times_K X)\cap H_a\subset
   \A^2_K\times_K \A^2_K=\A^4_K$, where $H_a$ is
the affine hyperplane $X_1-aX_2=0$ in $\A^4_K$.  In a neighborhood of the origin, $V_a$ is a curve, having the origin as a smooth $K$-point. Let $X_a$ be the unique irreducible
component of $V_a$ containing the origin.  Then $X_a(K)$ is infinite because $K$ is large and $X_a(K)$ is non-empty.
Thus there exists $(x_1,y_1,x_2,y_2)\in X_a(K)$ with $x_2\ne 0$.
Equivalently there exist $(x_1,y_1),(x_2,y_2)\in X(K)$ with $x_2\neq0$ and $a=x_1/x_2$, proving the claim. \qed

\medskip

\noindent{\bf Theorem~3.4.} {\it The function field $K$ of a smooth projective curve over a large field $F$ is quasi-free, of rank equal to the cardinality of $F$.}

\medskip

\noindent{\it Proof.}  In [HS, Theorem~4.3], it was shown that if $F$ is a very large field of cardinality $m$, and $K$ is the function field of a smooth projective $F$-curve, then every non-trivial finite split embedding problem for $K$ has $m$ proper regular solutions.  Hence the set of {\it all} proper solutions also has cardinality $m$; and thus $K$ is quasi-free [HS, Corollary~4.4].  
The result now follows from Proposition~3.3 above. \qed

\medskip

\noindent{\it Remark 3.5.} (a) As the proof of Theorem~3.4 shows, under the hypotheses of the theorem, every non-trivial finite split embedding problem for $K$ has exactly $m$ proper {\it regular} solutions, where $m = \card F$.  So this theorem strengthens Pop's result ([Po1], [Po2]) that if $K$ is the function field of a smooth projective curve over a large field, then every finite split embedding problem for $K$ has at least one proper regular solution.

(b) The property of being large (or PAC) can be regarded as complementary to the property of being Hilbertian (see [La], [FJ]).  Namely, consider a Galois branched cover $\phi:Y \to X = \A^1_K$.  If $K$ is Hilbertian, then there are infinitely many $K$-points of $X$ that remain prime in $Y$.  Meanwhile, to say that $K$ is PAC or large is to say that there are infinitely many $K$-points of $X$ that are totally split in $Y$ (in the latter case, assuming there is one such point).  Moreover, as for large fields, this property for curves implies a corresponding property in higher dimensions.  (Note also these properties are analogous to the two extremes in the Tchebotarev Density Theorem.)

(c) Remark (b) suggests introducing a notion of {\it very Hilbertian}; i.e.\ that for $Y \to X$ as in (b), the cardinality of the set of $K$-points of $X$ that remain prime in $Y$ is equal to the cardinality of $K$.  And in fact, the strategy of the proof of Theorem~3.3 also shows that every Hilbertian field is very Hilbertian.  Namely, if $\phi$ is generically given by a polynomial $f(x,y) \in K[x,y]$, consider for each $a \in K$ the variety $V_a \subset \A^4_K$ as in the proof of Theorem~3.3.  Then there exist $(x_1,x_2) \in \A^2(K)$ which remains prime in $V_a$, such that $x_2 \ne 0$; i.e.\ such that $f(x_1,Y)$ and $f(x_2,Y)$ are irreducible in $K[Y]$, with $x_1/x_2=a$.  The property of being very Hilbertian then follows.  

(d) If $K$ is a Hilbertian large field, then every finite split embedding problem over $K$ has a proper solution [Po2, Main Theorem~B] (since every finite split embedding problem over the function field of the $K$-line has a proper regular solution).  In fact, each such non-trivial embedding problem has infinitely many solutions (e.g.\ by replacing the kernel $N$ of the embedding problem by $N^n$ for arbitrarily large integers $n$, and then taking quotients of the solutions).  Since the properties of large and Hilbertian imply the properties of being very large and very Hilbertian, this suggests that a large Hilbertian field $K$ is quasi-free (and of rank equal to the cardinality of $K$).  Surprisingly, this is false, by an example of Jarden.  Namely, according to  Examples~3.1 and 3.2 of [Ja], there is a profinite group $G$ of uncountable rank that is projective and $\omega$-free but not free, and which is the absolute Galois group of a Hilbertian PAC (and hence large) field $K$.  Since $G = G_K$ is projective but not free, it cannot be quasi-free.   

(e) By another example (also due to Jarden), it is also possible for a large Hilbertian field $K$ to be quasi-free of rank strictly smaller than the cardinality of $K$.  Namely, by [FJ, Theorem~23.1.1], there is an uncountable field $K$ which is PAC and whose absolute Galois group $G_K$ is free of countable rank.  So $K$ is large, and $G_K$ is quasi-free of countable rank (and in particular $\omega$-free).  Also $K$ is Hilbertian by Roquette's theorem [FJ, Theorem~27.3.3] because it is $\omega$-free and PAC.  So $K$ is as claimed.  
Note that combining this example with Remark (c) above exposes a subtle point: for such a Hilbertian field $K$ and any finite Galois extension $L$ of $K(x)$, there will be $\card K$ elements of $K$ for which the specialization of $L$ is irreducible; but the corresponding Galois field extensions of $K$ are not linearly disjoint (and {\it up to isomorphism} there are fewer than $\card K$ of them).  

\bigskip

\noindent {\bf Section 4. Main results.}

\medskip

This section contains the main results of this paper, viz.\ the freeness of the absolute Galois groups of the function field of a real curve without real points, of the maximal abelian extension of $\C((x,y))$, and of the maximal abelian extension of the function field of a curve over a finite field.  Each of these is stated in somewhat stronger form below.

\medskip

\noindent{\bf Theorem 4.1.} (``Geometric Shafarevich Conjecture'') {\sl Let $p$ be a prime and let $k$ be a subfield of $\bar \F_p$ (e.g.\ a finite field).  Let $F$ be a one-variable function field over $k$, and let $F\ab$ be its maximal abelian extension.  Then the absolute Galois group of $F\ab$ is free of countable rank.}

\medskip

\noindent{\it Proof.}  Let $\tilde F$ be the compositum of $F$ and $\bar \F_p$ in an algebraic closure of $F$.  Then $\tilde F$ is the function field of a smooth projective curve over $\bar \F_p$.  Moreover we have containments $F \subset \tilde F \subset F\ab \subset \tilde F\ab$; i.e.\ $F\ab$ is abelian over $\tilde F$.  By [Ha1] or [Po1], the absolute Galois group of $\tilde F$ is a free profinite group of countably infinite rank.  So the same holds for $F\ab$, by Proposition~3.2(b) above.  \qed 

\medskip

Recall that a field $K$ with algebraic closure $\bar K$ is {\it formally real} if $-1$ is not a sum of squares in $K$; and $K$ is {\it real closed} if it is a maximal element in the set of formally real subfields of $\bar K$.  
If $K$ is real closed then $K[\sqrt{-1}]$ is algebraically closed; and so the absolute Galois group of a real closed field is cyclic of order $2$.  According to [CT, p.360], a field is large if its absolute Galois group is a pro-$p$ group for some prime $p$; in particular, real closed fields are large.  (More generally, according to [Po2, pp.~18-19], ``pseudo-real closed'' fields are large because they satisfy a universal local-global principle.)

\medskip

\noindent{\bf Theorem 4.2.} {\sl Let $X$ be a smooth projective curve $X$ over a real closed field $R$ (e.g.\ $R = \R$), and let $K$ be the function field of $X$.  Then the absolute Galois group of $K$ is free if and only if $X$ has no $R$-points; and if it is free, its rank is equal to the cardinality of $R$.}

\medskip

\noindent{\it Proof.}  As noted above, every real closed field is large.  So Theorem~3.3 says that $K$ is quasi-free of rank equal to $m := \card R$.  Thus $K$ is free (necessarily of rank $m$) if and only if it is projective, by [HS, Theorem~2.1]. 

In general, the function field of an integral variety of dimension $d$ over a real closed field $R$ with no $R$-points has cohomological dimension $d$ [CP, Proposition~1.2.1]. So in our situation, if $X(R)$ is empty then $K$ has cohomological dimension $1$, which implies that it is projective (as noted at the beginning of Section~2 above).  

Conversely, if $X$ has an $R$-point, then it is classical that $K$ is not projective.  This can also be seen directly as follows:  
Let $C=R[\sqrt{-1}]$. 
If $K=R(X)$ is projective, then the $\Z/2$-Galois extension $C(X)/R(X)$ can be embedded in a $\Z/4$-Galois field extension $L/R(X)$ (since the kernel of $\Z/4 \to \Z/2$ is Frattini).  Since $R$ is large, $X(R)$ is infinite; so some $P \in X(R)$ is unramified in this extension.  But a decomposition group over $P$ would then surject onto $Z/2$ and thus be $\Z/4$, which is impossible since $\Z/4$ is not a Galois group over $R$. \qed 

\medskip

\noindent{\it Remarks 4.3.}  a) As an example of the theorem, the fraction field of $\R[x,y]/(x^2+y^2+1)$ has free absolute Galois group, of rank equal to the cardinality of $\R$.

b) The proof of [CP, Proposition~1.2.1] is due to Ax and relies on a result of Serre.  But in the proof above, only the dimension $1$ case of [CP, Proposition~1.2.1] is needed; and that case is more classical, essentially going back to Witt [Wi].

c) The above result suggests asking about the structure of the absolute Galois group $G_K$ of the function field of a real curve $X$ that has a real point.  According to [SS, Theorem~2.2], the locus of real curves $X$ of genus $g$ is a union of connected (in fact irreducible) real analytic subspaces $M_\R^{(g,k,\varepsilon)}$, where $k, \varepsilon$ are non-negative integers.  Moreover $0 \le k \le g+1$ and $0 \le \varepsilon \le 1$; and $M_\R^{(g,k,\varepsilon)}$ is non-empty if and only if either $\varepsilon = 0$ and $k \ne g+1$ or else $\varepsilon = 1$ and $k \equiv g+1 \, ({\rm mod}\ 2)$ [SS, Theorem~1.1].  Here $(g,k,\varepsilon)$ is called the {\it type} of the real curve $X$ of genus $g$, where $k$ is the number of connected real components (ovals) of $X(\R)$ and where $\varepsilon=0$ if and only if $X(\R)$ is connected in the real topology.  Actually, the number $c$ of connected components of $X(\R)$ is at most $2$ [Mi, p.46], so $\varepsilon = c-1$.  As a real curve $X$ varies in moduli with the type remaining constant, the tower of branched covers of $X$ can also be deformed (preserving the numbers $r_1$ and $2r_2$ of real and complex branch points), since the topology remains constant.  So the algebraic fundamental group of the corresponding punctured curve also remains constant, for a given type and given values of $r_i$.  These fundamental groups fit together compatibly, as $r_1, r_2$ vary.  Thus the absolute Galois group $G_K$ depends (up to isomorphism) only on the type of $X$; and we can denote this group by $G_{g,k,\varepsilon}$.  The above theorem shows that $G_{g,k,\varepsilon}$ is free if and only if $k=0$.  It would be interesting to know the structure of the profinite group $G_{g,k,\varepsilon}$ for $k>0$, and whether it depends just on $k$
or also on $g$ and $\varepsilon$.

\medskip

Finally, we turn to consideration of the absolute Galois group of
the maximal abelian extension $K$ of $k((x,y))$, where $k$ is a separably closed field of arbitrary characteristic.  As in the previous theorem, we prove that this is free by using that it is projective and quasi-free.  For projectivity, we rely on the following result:

\medskip

\noindent {\bf Theorem~4.4.} [COP] {\it Let $A$ be an excellent henselian two-dimensional local domain, with fraction field $K$
and separably closed residue field $k$, of equal characteristic $p \ge 0$.  Then the maximal abelian extension $K\ab$ of $K$ has cohomological dimension at most $1$, as does the maximal pro-prime-to-$p$ abelian extension $K'$ of $K$.}

\medskip

\noindent {\it Proof.}  The case $p=0$ was shown in Theorem~2.3 of [COP] (Theorem~2.2 in the preprint).  For $p>0$, we modify that proof (following a sketch provided by R.~Parimala):

A field of characteristic $p\ne 0$ has $p$-cohomological dimension at most 1 [Se, II 2.2 Proposition 3].  So it suffices to show that $\cd_\ell \le 1$ for all $\ell \ne p$.  Regard $K \subset K' \subset K\ab \subset \bar K$, where $\bar K$ is a separable closure of $K$.  Since the extension $K\ab/K'$ is algebraic,  
$\cd_\ell(K\ab) \le \cd_\ell(K')$ [Se, II 4.1 Proposition 10].  So it suffices to consider just the case of $K'$.
By [Se, II 2.3 Proposition 4] and [Se, II 1.2 Proposition 1], $\cd_\ell(K') \le 1$ for $\ell \ne p$ if and only if every finite separable extension $F/K'$ satisfies $\Br(F)(\ell)=0$, where $\Br(F)(\ell)$ denotes the $\ell$-primary part of the Brauer group $\Br(F) := H^2(F,\G_m)$.  So it suffices to show that for every finite separable extension $F/K'$ (contained in $\bar K$), every central simple $F$-algebra of exponent prime to $p$ is split.

{\it Case I: $F$ is Galois over $K$.}  An $F$-algebra as above is induced via base change from a central simple algebra $D/L$, of exponent $n$ prime to $p$, where $L/K$ is a finite, Galois field extension of $K$.  Let $d=[L:K]$ and write $d=d'p^m$ with $d'$ prime to $p$ and $m \ge 0$.  Thus $p$ does not divide $N := nd'$.  The henselian ring $A$ (and hence $K$, $L$, $F$) contains all prime-to-$p$ roots of unity since the residue field $k$ is separably closed of characteristic $p$.  

Let $B$ be the integral closure of $A$ in the field $L$; this is an excellent henselian two-dimensional local domain, whose residue field $\tilde k$ is purely inseparable over $k$ and hence is separably closed.  As in the proof of [COP, Theorem 2.2] we obtain: 

-- a Weil divisor $\Delta$ on $\Spec(B)$ that is invariant under $G :=
\Gal(L/K)$, containing all singular points of $\Spec(B)$ and all points of codimension 1 where the algebra $D$ ramifies;

-- a proper integral regular model $\pi : X \to \Spec(B)$, where $\pi$ is projective and birational, $L$ is the function field of $X$, and the reduced divisor $\pi^{-1}(\Delta)_{\red}$ on $X$ is a $G$-invariant divisor with normal crossings, of
the form $C + E$, where $C$ and $E$ are regular closed curves; 

-- a finite, $G$-invariant set $S$ of closed points of $X$ including all points of intersection of $C$ and $E$ and at least one point of each component of $C + E$;

-- a function $g \in L^*$ such that $\div_X(g) =
C+E+J$ where $J$ is a divisor whose support does not contain any point of $S$, whose norm $f = N_{L/K}(g) \in K^*$ has divisor given by
$\div_X(f) = d \cdot (C+E) + \sum_{\sigma \in G} \sigma J$.

Let $M = L(f^{1/N}) \subset \bar K$, and let $D_M$ be the extension of $D$ to $M$.  Since $L$ contains a primitive $N$th root of unity, $M$ is a cyclic extension of $L$.  Also $M \subset F$, since $L \subset F$ (by definition of $L$) and since $h := f^{1/N} \in K' \subset F$.
So in order to show that the given central simple algebra is split over $F$ it suffices to show that $D_M$ is split.  

Let $B_1$ be the integral closure of $B$ in $M$; and let $Y \to \Spec(B_1)$ be a regular integral proper model, equipped with a projection map $q : Y \to X$ compatible with $\Spec(B_1) \to \Spec(B)$.  The unramified Brauer group of $M$, consisting of classes that are unramified with respect to all discrete valuations of $M$, is contained in the Brauer group $\Br(Y) := H^2_\et(Y,\G_m)$ of $Y$, since $Y$ is a regular surface with function field $M$ [COP, Corollary 1.9].  So it suffices to show that $D_M$ is unramified at every codimension 1 point $y$ on $Y$.  As in [COP] there are several cases.  If $x := q(y)$ does not belong to $C + E$, or if $x$ is of codimension 2 on $X$ and is not an intersection point of $C$ and $E$, then this property follows as in the proof of [COP, Theorem~2.3].  If $x$ has codimension 1 on $X$ and belongs to $C + E$, then from $h^N = f \in M^*$ we obtain
$nd' \cdot \div_Y(h) = N \cdot \div_Y(h) = \div_Y(f) = d'p^m \cdot q^{-1}(C + E) + q^{-1}(\sum_{\sigma \in G} \sigma J)$; 
so $n$ divides the ramification index of $y$ over $x$ (using that $n$ is prime to $p$) and $D_M$ is unramified at $y$ (as in [COP]).
Finally, in the case that $x$ is a intersection point of $C \cap E \subset X$, we have $f = u\pi^d\delta^d \in {\cal O}_{X,x}$, where $u \in {\cal O}^*_{X,x}$ and where $\pi, \delta \in{\cal O}_{X,x}$ form a regular system of parameters respectively defining $C$ and $E$ locally.  As in [COP] it suffices to show that the
symbol $(\pi, \delta)_n$ vanishes when viewed as an element of $\Br(M_y)$, where $M_y$ is the fraction field of the henselization ${\cal O}^h_{Y,y}$.  For this it suffices to show that $p^m\cdot (\pi, \delta)_n = 0$ in $\Br(M_y)$, because $n$ is relatively prime to $p$.  Since units in the multiplicative
group of ${\cal O}^h_{X,x}$ are divisible by integers that are prime to $p$, we have $h^{nd'} = f = v^{nd'} \pi^{d'p^m} \delta^{d'p^m} \in {\cal O}^h_{Y,y}$ for some $v \in ({\cal O}^h_{X,x})^*$.
The residue field of ${\cal O}^h_{Y,y}$ contains the separably closed field $\tilde k$, and so the group of roots of unity in ${\cal O}^h_{Y,y}$ is $d'$-divisible.
Thus $(\pi\delta)^{p^m} = \rho^n$ for some $\rho \in M_y$. So in $_n\Br(M_y) = H^2_{\et}(M_y, \mu_n) \simeq H^2_{\et}(M_y, \mu_n^{\otimes 2})$, we obtain as desired
$p^m\cdot (\pi, \delta)_n = p^m \cdot (\pi, \pi^{-1})_n + (\pi, \rho^n)_n = 0 + 0 = 0$.

{\it Case II: General case.} Let $M \subset \bar K$ be the Galois closure of $F$ over $K$; this is finite over $F$.  By Case I, $\Br(M)(\ell) = 0$; so $[M:F]\Br(F)(\ell) = 0$.  Choosing an isomorphism of $\Z[1/\ell]/\Z$
with the $\ell$-power roots of unity of $F$, the Merkuriev-Suslin theorem [MS] gives an isomorphism 
$K_2(E) \otimes (\Z[1/\ell]/\Z) \approx \Br(E)(\ell)$; so $\Br(E)(\ell)$ is $\ell$-divisible.  But being a pro-$\ell$-group, $\Br(E)(\ell)$ is also $r$-divisible for every integer $r$ that is prime to $\ell$.  So $\Br(E)(\ell)$ is divisible; and hence is trivial, being $[M:F]$-torsion. \qed

\medskip

\noindent{\it Remark.}  The above proof breaks down in the unequal characteristic case, where $\char\, K = 0$ and $\char\, k = p \ne 0$, because of the need in that case to show that $\cd_p \le 1$.

\medskip

Using the above result, we obtain:

\medskip

\noindent{\bf Theorem~4.5.} {\sl Let $k$ be a field and let $K$ be the maximal abelian extension of $k((x,y))$, with absolute Galois group $G_K$.  

a) Then $G_K$ is quasi-free of rank equal to the cardinality of $K$.

b) If $k$ is separably closed, then the absolute Galois group of $K$ is a free profinite group of rank equal to the cardinality  of $K$.}

\medskip

\noindent{\it Proof.}  a) According to [HS, Theorem~5.1], the absolute Galois group of $K_0 := k((x,y))$ is quasi-free of rank equal to $\card K_0$ (even without any assumptions on $k$).  By Proposition~3.1(b), it follows that the absolute Galois group of $K = K_0\ab$ is also quasi-free of rank $\card K_0$.  But $K$ and $K_0$ have the same cardinality; so the assertion follows.

b) By [HS, Theorem~2.1], a profinite group is free of infinite rank $m$ if and only if it is projective and is quasi-free of that rank.  
As noted before, $G_K$ is projective if and only if $K$ has cohomological dimension $1$; and that latter property holds by Theorem~4.4.  So the assertion follows from part (a).  \qed

\medskip

\noindent{\it Remark.}  The above proof of Theorem 4.5(b) relies on 4.5(a), and hence on Proposition~3.1 (and thus Theorem~2.4).  But if one is willing to bypass 4.5(a), one can prove a weaker version of 4.5(b) --- that $K$ is $\omega$-free --- without relying on those other results.  Namely, one can proceed as follows, as suggested by M.~Jarden:   By a theorem of Weissauer [FJ, Theorem~15.4.6], $K_0 = k((x,y))$ is Hilbertian, being the fraction field of the two dimensional Krull domain $k[[x,y]]$.  So its maximal abelian extension $K$ is also Hilbertian, by a theorem of Kuyk [FJ, Theorem~16.11.3].  Thus every finite split embedding problem for $K$ with an abelian kernel has a proper solution, by a theorem of Ikeda [FJ, Proposition~16.4.5].  Since $G_K$ is projective by Theorem~4.4 above (using that $k$ is separably closed), every finite embedding problem for $K$ is dominated by a finite split embedding problem; and so solving any finite embedding problem for $K$ can be reduced to solving a finite sequence of finite split embedding problems each of which has a minimal normal subgroup as its kernel.  So it suffices to show that such embedding problems have proper solutions.  If the kernel of such an embedding problem is abelian, then we are done by the theorem of Ikeda cited above.  Otherwise, the kernel of the embedding problem is a product of finitely many isomorphic non-abelian finite simple groups [As, Chap.~3, 8.3, 8.2].  This embedding problem for $K$ is induced by a finite split embedding problem for some finite extension $K_1$ of $K_0$ that is contained in $K$.  But $K_0$ is quasi-free by [HS, Theorem~5.1]; and hence so is $K_1$, by [RSZ].  So there is a proper solution to the embedding problem for $K_1$; and this induces a proper solution to the embedding problem over $K$ because of linear disjointness, since $K$ is abelian over $K_1$ whereas the kernel of the embedding problem has no non-trivial abelian quotients.  So $K$ is $\omega$-free.

\bigskip

\noindent{\bf References.}

\medskip

\noindent[As] M.~Aschbacher.  ``Finite Group Theory'', 2nd ed.  Cambridge Univ.\ Press, 2000.

\smallskip

\noindent[Bo] F.~Bogomolov.  On the structure of Galois groups of the fields of rational function.  Proc.\ Symp.\ Pure Math., {\bf 58.2} (1995), 83-88.

\smallskip

\noindent[CT] J.-L.~Colliot-Th\'el\`ene.  Rational connectedness and Galois covers of the projective line.  Annals of Math., {\bf 151} (2000), 359-373.

\smallskip

\noindent[COP] J.-L.~Colliot-Th\'el\`ene, M.~Ojanguren, and R.~Parimala.
Quadratic forms over fraction fields of two-dimensional Henselian rings and Brauer groups of related schemes.  In: ``Proceedings of the International Colloquium on Algebra, Arithmetic and Geometry'', Tata Inst.\ Fund.\ Res.\ Stud.\ Math., vol.~16, pp.185-217, Narosa Publ.\ Co., 2002.   Also available as preprint at $\langle$http://www.math.u-psud.fr/$\sim$colliot/CTOjPa22may01.ps$\rangle$.

\smallskip

\noindent[CP]  J.-L.~Colliot-Th\'el\`ene and R.~Parimala.  Real components of algebraic varieties and \'etale cohomology.  Invent.\ Math., {\bf 101} (1990), 81-99.

\smallskip

\noindent[FJ] M.~Fried and M.~Jarden.  ``Field Arithmetic'', 2nd ed.  Ergebnisse der Math., vol.~11.  Springer, 2005.

\smallskip

\noindent[HJ] D.~Haran, M.~Jarden.  Regular split embedding problems over complete valued fields.  Forum Mathematicum {\bf 10} (1998), 329-351.

\smallskip

\noindent[Ha1] D.~Harbater. Fundamental groups and embedding problems in characteristic $p$.  In: ``Recent Developments in the Inverse Galois Problem'' (M. Fried, et al., eds.), AMS Contemp. Math. Series, vol.~186, 1995, pp.~353-369.
\smallskip

\noindent[Ha2] D.~Harbater.  Shafarevich conjecture.  In ``Encyclopaedia of Mathematics'', Supplement III.
Managing Editor: M. Hazewinkel, Kluwer Academic Publishers, 2002,
pp.360-361.

\smallskip

\noindent[Ha3] D.~Harbater.  Patching and Galois theory.  In ``Galois Groups and Fundamental Groups'' (L. Schneps, ed.), MSRI Publications series, vol.41, Cambridge Univ. Press, 2003, pp.313-424.

\smallskip

\noindent[HS] D.~Harbater and K.~Stevenson.
Local Galois theory in dimension two. Advances in 
Math.\ (special issue in honor of M. Artin's 70th birthday), 
{\bf 198} (2005), 623-653.

\smallskip

\noindent[Iw] K.~Iwasawa.  On solvable extensions of algebraic number fields.
Annals of Math.\ {\bf 58} (1953), 548-572.

\smallskip

\noindent[Ja] M.~Jarden. On free profinite groups of uncountable rank.  In ``Recent developments in the inverse
Galois problem'' (M.~Fried, et al., eds.), AMS Contemp.\ Math. Series, vol.~186, 1995, pp.~371-383.

\smallskip

\noindent[La] S.~Lang, ``Fundamentals of Diophantine Geometry'', Springer, 1983.

\smallskip

\noindent[Mi] G.~Mikhalin.  Adjunction inequality for real algebraic curves.  Math.\ Res.\ Letters, {\bf 4} (1997), 45-52.

\smallskip

\noindent[MS] A.S.~Merkuriev and A.A.~Suslin.
$K$-cohomology of Severi-Brauer varieties and the norm residue homomorphism (in Russian).
Izv.\ Akad.\ Nauk SSSR Ser.\ Mat., {\bf 46} (1982), 1011-1046, 1135-1136.  English translation: Math.\ USSR-Izv.\ {\bf 21} (1983),  307-340.

\smallskip

\noindent[Po1] F.~Pop.  \'Etale Galois covers of affine smooth curves.  Invent.\ Math., {\bf 120} (1995), 555-578.

\smallskip

\noindent[Po2] F.~Pop.  Embedding problems over large fields.  Ann.\ Math., {\bf 144} (1996), 1-34.

\smallskip

\noindent[RSZ] L.~Ribes, K.~Stevenson, and P.~Zalesskii.  On quasifree profinite groups.  2006 manuscript.  To appear in Proceedings of the AMS.

\smallskip

\noindent[SS] M.~Sepp\"al\"a and R.~Silhol.  Moduli spaces for real algebraic curves and real abelian varieties.  Math.\ Zeitschrift {\bf 201}  (1989), 151-165.

\smallskip

\noindent[Se] J.-P.~Serre. ``Cohomologie Galoisienne'', 4th ed.,
Lec.\ Notes in Math., {\bf 5}, Springer, 1973.

\smallskip

\noindent[Wi] E.~Witt.  Zerlegung reeller algebraischer Funktionen in Quadrate, Schiefk\"orper \"uber reellem Funktionenk\"orper.  J.\ Reine Angew.\ Math.\ {\bf 171}, 4-11 (1934).

\end